# Detection of Approaching Critical Transitions in Natural Systems Driven by Red Noise

Andreas Morr* and Niklas Boers

*Earth System Modelling, School of Engineering and Design, Technical University Munich, Germany*
*and Research Domain IV-Complexity Science, Potsdam Institute for Climate Impact Research, Germany*

Detection of critical slowing down (CSD) is the dominant avenue for anticipating critical transitions from noisy time-series data. Most commonly, changes in variance and lag-1 autocorrelation [AC(1)] are used as CSD indicators. However, these indicators will only produce reliable results if the noise driving the system is white and stationary. In the more realistic case of time-correlated red noise, increasing (decreasing) the correlation of the noise will lead to spurious (masked) alarms for both variance and AC(1). Here, we propose two new methods that can discriminate true CSD from possible changes in the driving noise characteristics. We focus on estimating changes in the linear restoring rate based on Langevin-type dynamics driven by either white or red noise. We assess the capacity of our new estimators to anticipate critical transitions and show that they perform significantly better than other existing methods both for continuous-time and discrete-time models. In addition to conceptual models, we apply our methods to climate model simulations of the termination of the African Humid Period. The estimations rule out spurious signals stemming from nonstationary noise characteristics and reveal a destabilization of the African climate system as the dynamical mechanism underlying this archetype of abrupt climate change in the past.

## I. INTRODUCTION

The phenomenon of critical slowing down (CSD), which occurs in the advent of critical transitions induced by certain bifurcations, is an essential observational characteristic in the analysis of dynamical systems and for the anticipation of such transitions [1]. If caused by the approach of a codimension-1 bifurcation, the vanishing of a stable equilibrium point and the resulting transition will be preceded by a gradual decline of the linearized restoring rate of said equilibrium. This decline, in turn, leads to a weaker and slower response to perturbations, i.e., higher variance and autocorrelation in time. This result can indeed be shown analytically for fold-type bifurcations driven by small, additive white noise with standard deviation $\sigma$ and the drift of the linearized dynamics denoted by $\lambda$: for the variance, $\langle x^2 \rangle = \sigma^2/(2\lambda)$, and for the autocorrelation, $AC(\tau) = \exp(-\lambda\tau)$. Insofar as the assumption that high-dimensional complex systems such as Earth system components are prone to bifurcation-induced tipping is justified, CSD is expected to occur in the dynamics leading up to these events [2–4]. This idea has spurred interest in the development of so-called CSD indicators or early warning signals (EWS), i.e., estimators of local system stability that allow one to anticipate bifurcation-induced transitions [3,5–8]. However, the applicability of such estimators will depend on whether the actual system's dynamics are approximated well by the simple low-dimensional model used to derive them. This question pertains both to the approximation of the deterministic equilibrium dynamics [9–12] and the representation of omitted dimensions via a stochastic component in terms of noise [13–15]. In the most reductive model for fold-type bifurcations, a one-dimensional observable $X$ of the system is assumed to remain close to equilibrium and thus experience approximately linear restoring forces,

$$dX_t^{(w)} = -\lambda X_t^{(w)} dt + dW_t.$$

The linear restoring rate $\lambda$ will then vanish gradually as the system approaches the critical forcing value of the fold bifurcation. Perturbations to the system are modeled as additive white noise $dW$, with $W$ being a Wiener process. Its use assumes temporal independence of the perturbations inflicted on the system by the unresolved dynamics. However, many physical systems exhibit memory effects or persistence in their unresolved dynamics. In particular, the Mori-Zwanzig formalism implies that if an effective stochastic dynamic equation of a high-dimensional system is derived as the projection to a low-dimensional space of observed variables, the interactions between resolved and unresolved variables lead to non-Markovian dynamics [16,17]. To represent the memory, a model driven by red noise,

$$dX_t = -\lambda X_t dt + U_t dt,$$

with an Ornstein-Uhlenbeck process $U$, is more suitable [18,19]. While other continuous-time noise models with positive correlation in time exist [20,21], the specific frequency characteristics of red noise make it the most appropriate for application to many physical systems, including the Earth's climate [22–24]. Numerous techniques exist for assessing system stability under the influence of white noise [6,25,26]. In contrast, the red noise case addressed here has so far only been approached from the standpoint of discrete-time models [13,27,28]. In the following, we see that the white noise case can be obtained as a parameter limit of the red noise case. We introduce two novel stability indicators designed to simultaneously suitable for the red and white noise cases and compare their performance to the well-established variance and lag-one autocorrelation in the general case of nonstationary time-correlated driving noise. We also discuss the applicability of two existing discrete-time methods developed for such nonstationary noise, presented in Refs. [13,27], and assess their performance in the continuous-time setup. Lastly, we apply the novel methods to time-series data of the abrupt transition ending the African Humid Period, which was recently reproduced in simulations with a global climate model [29].

*andreas.morr@tum.de

## II. METHODS

### A. Linearly restoring process under red noise forcing

We first linearize the dynamics of the observable $x$ around a fixed point $x^*$:

$$\dot{x}(t) = f(x; a) \approx \partial_x f(x^*(a); a)(x(t) - x^*)$$
$$=: -\lambda(a)(x(t) - x^*).$$

The dynamics described via $f$ are mutable through the external parameter $a$ and are assumed to be autonomous. If the global dynamics $f(x; a)$ are that of a generic fold bifurcation located at a certain value $\bar{a}$, then the linearized restoring rate $\lambda(a) > 0$ of the initial state will decrease and eventually vanish: $\lambda(\bar{a}) = 0$.

The following formulas give a general model of a system driven by positively correlated noise, which is of particular interest when considering CSD in physical systems with the corresponding dynamics,

$$dX_t = -\lambda X_t dt + \kappa U_t dt, \qquad X_0 = 0, \qquad (1a)$$

$$dU_t = -\theta U_t dt + dW_t, \qquad U_0 = 0, \qquad (1b)$$

where $W$ is a Wiener process on the filtered probability space $(\Omega, \mathcal{F}, (\mathcal{F}_t)_{t \in \mathbb{R}}, \mathbb{P})$. For comparison, we also consider the model forced by white noise:

$$dX_t^{(w)} = -\lambda X_t^{(w)} dt + \sigma dW_t, \qquad X_0^{(w)} = 0. \qquad (2)$$

The solutions to the stochastic differential equations (SDEs) (1) and (2) are

$$X_t = \kappa \int_0^t \exp(-\lambda(t - s)) U_s ds, \qquad (3)$$

$$X_t^{(w)} = \sigma \int_0^t \exp(-\lambda(t - s)) dW_s, \qquad (4)$$

respectively. Note that all parameters of the model are *a priori* assumed to be constant because we are interested in the stationary characteristics of the observable at any given distance from the bifurcation point. From these characteristics, we derive suitable estimators of the constant linear restoring rate $\lambda$.

Both Eqs. (3) and (4) are asymptotically stationary Gaussian processes. There exist initial distributions for $X_0^{(w)}$ and $(X_0, U_0)$, respectively, such that they are stationary for all $t \geq 0$. The stationary characteristics of $X^{(w)}$ are well known since it is itself an Ornstein-Uhlenbeck process. They are given in Table I. For the red-noise-driven process $X$, we derive these characteristics via the corresponding Lyapunov equation (see Table II). We include explicit calculations in the Supplemental Material (SM) [30].

We further observe that the stochastics of the white-noise-driven process $X^{(w)}$ are the limits in the distribution of the stochastics of the red-noise-driven process $X$ in the case that $\kappa, \theta \to \infty$ and $(\kappa/\theta) \to \sigma$. This specific convergence is an example of a more general model convergence, which has been discussed extensively in the literature [31]. It will later be of practical use when employing estimators for system stability that are sensitive to such a limit. We will henceforth only consider $X$ as the general model and implicitly include the setting $X^{(w)}$ as a limit case.

An important characteristic of the stationary distribution of $X$ is that all quantities are symmetric with respect to a swapping of $\lambda \leftrightarrow \theta$. This characteristic will be particularly relevant when trying to infer information about one of the two in isolation. Analyzing the behavior of the quantities introduced in Table II under changes in the three parameters $\lambda$, $\theta$, and $\kappa$, the risk of spurious CSD indications, and hence false alarms, becomes evident (see Table III). On the other hand, it is easy to imagine that simultaneous trends in the parameters could cause the respective observable quantity to remain constant, leading to missed alarms. We will later refer to this second case as a masking of CSD.

TABLE I. Stationary characteristics of the white-noise-driven process $X^{(w)}$ defined by Eq. (2). The continuous-time (c.t.) and discrete-time (d.t.) power spectral densities (PSDs) are defined as the time-scaled c.t. and d.t. Fourier transforms taken in their squared amplitude, respectively (see SM [30] for further discussion).

| | |
|---|---|
| Variance | $\sigma^2/(2\lambda)$ |
| AC($\tau$) | $\exp(-\lambda|\tau|)$ |
| c.t. PSD $S(\omega)$ | $\sigma^2/(\lambda^2 + \omega^2)$ |
| d.t. PSD $S^{(1)}(\omega)$ | $\frac{\sigma^2}{2\lambda} \frac{\sinh(\lambda)}{\cosh(\lambda) - \cos(\omega)}$ |

TABLE II. Stationary characteristics of the red-noise-driven process $X$ defined by Eq. (1). All quantities associated with the process $X$ are symmetric with respect to a swapping of $\lambda \leftrightarrow \theta$. The quantities calculated for the case $\lambda = \theta$ coincide with the well-defined limit $\theta \to \lambda$ of the case $\lambda \neq \theta$.

| | $\lambda \neq \theta$ | $\lambda = \theta$ |
|---|---|---|
| Variance | $\kappa^2/(2\lambda\theta(\lambda+\theta))$ | $\kappa^2/(4\lambda^3)$ |
| AC($\tau$) | $(\lambda\exp(-\theta|\tau|) - \theta\exp(-\lambda|\tau|))/(\lambda-\theta)$ | $(1+\lambda|\tau|)\exp(-\lambda|\tau|)$ |
| c.t. PSD $S(\omega)$ | $\kappa^2/((\theta^2+\omega^2)(\lambda^2+\omega^2))$ | $\kappa^2/(\lambda^2+\omega^2)^2$ |
| d.t. PSD $S^{(1)}(\omega)$ | $\frac{\kappa^2}{2\lambda\theta(\lambda^2-\theta^2)}\left(\frac{\lambda\sinh(\theta)}{\cosh(\theta)-\cos(\omega)} - \frac{\theta\sinh(\lambda)}{\cosh(\lambda)-\cos(\omega)}\right)$ | $\frac{\kappa^2}{4\lambda^3}\frac{\cosh(\lambda)(\lambda\cos(\omega)+\sinh(\lambda))-\sinh(\lambda)\cos(\omega)-\lambda}{(\cosh(\lambda)-\cos(\omega))^2}$ |

We note that there exists an ARMA(2,1) representation of the process $X$:

$$X_{k+1} \overset{d}{=} \left(\exp(-\lambda) + \exp(-\theta)\right)X_k - \exp\left(-(\lambda+\theta)\right)X_{k-1} + \sigma_0 z_k + \sigma_1 z_{k-1}, \tag{5}$$

where the $z_k$ are independent identically distributed (i.i.d.) unit normal random variables and the constants $\sigma_0$ and $\sigma_1$ are unwieldy yet may be explicitly computed by solving the appropriate system of correlation equations. A CSD indicator relying on the ARMA($p$, $q$) best model fit to data with no specific a priori fixed model structure has recently been proposed [28]. The above considerations on the red-noise-driven process $X$ imply that this method should be sensitive to CSD in this model. At the same time, the symmetry in the parameters with respect to $\lambda \leftrightarrow \theta$ implies a risk of spurious indications in the case of nonstationary noise, much like the conventional methods of variance and AC(1). Therefore, we will not include this approach in our later comparisons of indicator performances.

TABLE III. Behavior of the quantities in Table II when one of the three parameters is taken to the respective limit in the first column. The true CSD to anticipate a critical transition is present in the $\lambda \to 0$ case. Considering the $\theta \to 0$ case, the potential for a spurious indication of CSD is evident, while the AC(1) will reveal the $\kappa \to \infty$ case as spurious. Spectral reddening refers to the value of the PSD at low frequencies, in this case, $\omega = 0$. The behavior is equivalent for the discrete-time and continuous-time PSD. All of the increases are strictly monotonic. The symbol "—" refers to the independence of the AC(1) from the parameter $\kappa$.

| | Variance | AC(1) | Spectral reddening |
|---|---|---|---|
| $\lambda \to 0$ | $\nearrow\infty$ | $\nearrow 1$ | $\nearrow\infty$ |
| $\theta \to 0$ | $\nearrow\infty$ | $\nearrow 1$ | $\nearrow\infty$ |
| $\kappa \to \infty$ | $\nearrow\infty$ | — | $\nearrow\infty$ |

### B. Estimators of system stability $\lambda$

Perhaps the most common indicators in use for the detection of CSD are increases in variance and AC(1) of the observable $X$ [7,32,33]. As we have established above, both quantities will monotonically increase in the event of a decreasing linear restoring rate $\lambda \to 0$ under either red or white noise forcing. We first present well-established estimators for these two quantities before introducing one known and two novel estimation techniques for inferring information about the linear restoring rate. We see that, for each of the estimators, the white noise limit $\kappa, \theta \to \infty$ and $(\kappa/\theta) \to \sigma$ is well defined and consistent with the quantities obtained when applying the techniques to the white noise model. In our setup, the white noise case is hence a special case of the more general red noise model. Even without an a priori model decision on whether the noise is white or red, the introduced estimators are generally applicable. We discuss the application to time-series samples with a dimensionless time step $\Delta t = 1$, though the methods are, in principle, applicable to time series with any constant time step. Proofs for the applicability of the conventional estimation methods as well as a comparison of the quality of the estimators in terms of their sample spread for different parameter settings can be found in the SM [30], as can a numerical analysis of their distributional convergence in a central limit theorem fashion. While we do not prove a corresponding result, the numerical results suggest an underlying convergence property of our new estimators.

#### 1. Variance

A consistent estimator for the variance is

$$\widehat{\mathrm{Var}}_N := \frac{1}{N}\sum_{k=0}^{N-1} X_k^2,$$

converging in probability to the quantity determined in the previous section,

$$\widehat{\mathrm{Var}}_N \overset{\mathbb{P}}{\underset{N\to\infty}{\longrightarrow}} \frac{\kappa^2}{2\lambda\theta(\lambda+\theta)}.$$

### 2. Lag-1 autocorrelation

A consistent estimator for the lag-$\tau$ autocorrelation for $\tau < N$ is

$$\widehat{AC(\tau)}_N := \frac{N}{N-\tau} \frac{\sum_{k=0}^{N-\tau-1} X_k X_{k+\tau}}{\sum_{k=0}^{N-1} X_k^2}, \quad (6)$$

also converging in probability as

$$\widehat{AC(\tau)}_N \xrightarrow[N\to\infty]{\mathbb{P}} \frac{\lambda \exp(-\theta|\tau|) - \theta \exp(-\lambda|\tau|)}{\lambda - \theta}.$$

### 3. Generalized least squares estimator

There exist three notable studies regarding the detection of CSD under the influence of nonstationary time-correlated noise [13,27], with the third requiring explicit external knowledge of the noise characteristics [34]. Boettner and Boers [13] and Boers [27] build on the discrete-time model of an AR(1) process, in turn, driven by an AR(1) process

$$Y_{k+1} = \varphi Y_k + c V_k,$$
$$V_{k+1} = \rho V_k + z_k,$$

where the $z_k$, $k \in \mathbb{N}$ are i.i.d. unit normal random variables. Here, an increase towards 1 of the autoregressive parameter $\varphi$ would be indicative of a destabilization of the underlying dynamics and hence a sign of CSD. Rearranging these discrete-time evolution equations, one arrives at the following ARMA(2,0) model for $Y$:

$$Y_{k+1} = (\varphi + \rho) Y_k - \varphi \rho Y_{k-1} + c z_{k-1}. \quad (7)$$

Recalling the ARMA(2,1) representation of the continuous-time process $X$ given in Eq. (5), it is clear that because $\sigma_1 \neq 0$, the marginal distributions of $X_t$ and $Y_k$ will differ in their moments and correlations. Nevertheless, it is conceivable that the methods developed for the discrete-time model might deliver satisfactory results even on data from the continuous-time case.

The unbiased estimator for $\varphi$ introduced by Boettner *et al.* in Ref. [13] does not appear to be applicable. The reason is that, even when applied to time-series data generated through the intended model, Eq. (7), the algebraic expression of the estimator is not well defined on a set of positive probability, only performing well on time series much longer than the ones considered here. Applying the method to data of the continuous-time process $X$ seems to exacerbate this issue, effectively making the interpretation of the estimator results impossible. Thus, we do not consider it for further analysis.

The method proposed and implemented by Boers in Ref. [27] builds on regressing observed increments on the left-hand side against the system state on the right-hand side:

$$Y_{k+1} - Y_k = (\varphi - 1) Y_k + c V_k.$$

Instead of an ordinary least squares model suitable for white noise, the AR(1) structure of $V$ is taken into account. To this effect, the PYTHON module statsmodels and its class GLSAR are used. The resulting estimate $\hat{\varphi}$ is taken as a stability estimator, and its increase is taken as a CSD indicator. Comparing the ARMA models (5) and (7), one could assume that the underlying value of $\varphi$ should be approximately $\exp(-\lambda)$. However, investigating the distribution of $\hat{\varphi}$, its mean seems to significantly differ from this value (see Fig. 3.1 in the SM [30]). This finding is again due to the different ARMA structures.

### 4. Fitting to the observed autocorrelation structure

The symmetry of the stationary distribution of $X$ with respect to exchanging $\lambda$ and $\theta$ implies that explicit information about the parameters cannot be inferred from one-dimensional time-series statistics of variance and AC(1). The first novel method we propose circumvents this problem by including multiple estimated moments in the assessment.

Estimating the autocorrelation structure (ACS) of the observed process $X$ via the already-established estimator in Eq. (6), we find a tuple $(\hat{\lambda}^{(ACS)}, \hat{\theta}^{(ACS)})$ that constitutes the best model fit in the sense that the mean squared error between the observed $\widehat{AC(\tau)}_N$ and the theoretically computed $AC_{\lambda,\theta}(\tau)$ corresponding to the red noise model with these parameters (see Table II) is minimized. This case can be realized numerically by running a minimization function on the mean squared error:

$$(\hat{\lambda}_N^{(ACS)}, \hat{\theta}_N^{(ACS)}) = \arg\min_{\substack{(\lambda,\theta) \in \mathbb{R}_+^2 \\ \lambda < \theta}} \sum_{\tau=1}^{\tau_{max}} \left( \widehat{AC(\tau)}_N - AC_{\lambda,\theta}(\tau) \right)^2.$$

Since the set of arguments $\{(\lambda, \theta) \in \mathbb{R}^2 | \lambda < \theta\}$ is an open set, the minimum of the above squared error does not exist *a priori*. In the numerical implementation, either a local minimum is found or the estimation attempt fails. To include the white noise limit, the edge case of $\theta = \infty$ should be caught during the optimization and interpreted appropriately (see also SM [30] for a more detailed description of the numerical routine).

Though the idea of performing parameter estimation through the method of moment fitting is not new [35], so far it has not been applied to this specific problem. Note that we have made the model assumption that the correlation time $1/\theta$ of the noise component is always shorter than the correlation time $1/\lambda$ induced by the (locally) linear restoring dynamics. While the method is also applicable without this assumption, we have to bear in mind that some outside knowledge about the relation of the two parameters is required in order to distinguish the trends observed in them. A relative timescale separation in the noise and the dynamics of interest is a common assumption in many applied fields such as climate science [22,36]. Furthermore, if the linear restoring rate $\lambda$ indeed undergoes a decrease towards zero, at some point, it will fall below the value of $\theta$.

Choosing a "good" maximum $\tau_{\max}$ of evaluated lags is not easy to achieve comprehensively. Estimations of the autocorrelation deteriorate with increasing $\tau$, and the exponential decay of the theoretical model ACS implies that, for high lags, the change in neighboring lags is negligible. For all applications we are considering, a choice of $\tau_{\max} = 3$ delivers satisfactory results. A proof for the convergence of this estimator could be obtained by adapting the proof of Lemma 3.4 in Ref. [35], though we do not attempt it here.

### 5. Fitting to the observed power spectral density

Similarly, determining the model with the least mean squared error between the theoretically computed model PSD and the observed PSD suggests a choice of $(\hat{\lambda}^{(\mathrm{PSD})}, \hat{\theta}^{(\mathrm{PSD})}, \hat{\kappa}^{(\mathrm{PSD})})$. In this case, the observed PSD is the squared absolute value of the discrete-time Fourier transform of the data:

$$\widehat{S^{(1)}(\omega)}_N = \left| \frac{1}{\sqrt{N}} \sum_{k=0}^{N-1} \exp(-i\omega k) X_k \right|^2 .$$

In contrast to the previous ACS case, the discrete-time PSD is not equal to the continuous-time PSD, and it is imperative to choose the former, given in Table II.

The discrete-time PSD will be a periodic function classically probed on the frequencies $F = \{(2\pi l / N) | l = 1, \ldots, N/2 - 1\}$ if $N$ is even and $F = \{(2\pi l / N) | l = 1, \ldots, (N-1)/2\}$ if $N$ is odd. In order to weight the entire frequency range more evenly, taking the logarithm of the observed and expected PSD is advantageous. Averaging over neighboring frequencies to smooth out the fitting target may also improve the quality of the estimations. The estimator is then given by

$$\left( \hat{\lambda}_N^{(\mathrm{PSD})}, \hat{\theta}_N^{(\mathrm{PSD})}, \hat{\kappa}_N^{(\mathrm{PSD})} \right)$$
$$= \underset{\substack{(\lambda, \theta, \kappa) \in \mathbb{R}_+^3 \\ \lambda < \theta}}{\arg\min} \sum_{\omega \in F} \left( \log\left( \widehat{S^{(1)}(\omega)}_N \right) - \log\left( S_{\lambda, \theta, \kappa}^{(1)}(\omega) \right) \right)^2 .$$

Much like in the formulation of the estimators relying on the ACS, the set of arguments is open, and possible infima of the squared error on the boundary should be interpreted correctly in implementations. This case is, again, particularly relevant for the detection of the white noise limit $\theta, \kappa \to \infty$ and $(\kappa/\theta) \to \sigma \in \mathbb{R}_+$ (see, again, SM [30]).

This method bears similarities with the ratio of spectra (ROSA) method recently proposed in Ref. [34]. In their approach, Clarke *et al.* also performed a least square error fit between two PSDs, but they relied on dividing out the PSD of the driving noise, which needs to be known *a priori*. Therefore, this method is not suited to infer information about the stability of the system from the observable alone. In their practical implementation, they reverted to the continuous-time PSD as a theoretical fitting target. This choice can lead to considerable biases due to the mismatched behavior of the discrete-time PSD for frequencies close to the Nyquist frequency.

Using the PSD instead of the ACS as a model-fit target has two practical advantages in our context. First, since we are using the entire relevant frequency domain, we are not faced with having to fix another degree of freedom in the estimation. In the ACS method, this was the number of included lags $\tau_{\max}$. Second, since only the omitted zero frequency entry of the PSD is sensitive to a shift of the time series by a constant, the method is considerably more stable with respect to prior centering and low-order detrending. This feature is particularly relevant in applications where the approximate Ornstein-Uhlenbeck residual first has to be separated from a slow deterministic trend.

## III. RESULTS

### A. Comparison of the indicators

To compare the performance of the proposed indicators, we first formulate a general application setting. This setting will describe the range of possible parameter evolutions we posit for some real-world case of detecting CSD. In the classical setting, we would assume the white noise limit $\kappa / \theta \to \sigma$ and further assume that $\sigma$ is fixed during the time of observation. In that case, none of the indicators is prone to spurious indication or masking of CSD. However, in the general noise case, not only the parameter of interest, i.e., $\lambda$, but also the noise parameters $\theta$ and $\kappa$ change in time. In this case, the two conventional indicators likely give false-positive (spurious) or false-negative (masking) results. To quantitatively compare these pitfalls, we perform a disjoint window analysis on data from a large range of randomly drawn parameter settings. We check the resulting series of estimations for a positive Kendall's $\tau$ in $\widehat{\mathrm{AC}}(1)$, $\widehat{\mathrm{Var}}$, $\hat{\varphi}$, $-\hat{\lambda}^{(\mathrm{ACS})}$, and $-\hat{\lambda}^{(\mathrm{PSD})}$, respectively, each suggesting CSD. We then plot the receiver-operating characteristic (ROC) for each indicator to compare their ability to discern the cases with a truly decreasing linear restoring rate $\lambda$ from those where no change takes place.

The general model setting in which we probe our estimators is defined by

$$\mathrm{d}X_t = -\lambda(t) X_t \mathrm{d}t + \kappa(t) U_t \mathrm{d}t, \qquad X_0 = 0, \quad (8\mathrm{a})$$

$$\mathrm{d}U_t = -\theta(t) U_t \mathrm{d}t + \mathrm{d}W_t, \qquad U_0 = 0. \quad (8\mathrm{b})$$

Linearizing the equilibrium dynamics via $\lambda$ is expected to be a good enough approximation to justify this setup, replacing an actual codimension-1 bifurcation. Since the parameters are now deterministic functions of time, our considerations about the formerly stationary process $X$ do not hold exactly anymore. Nevertheless, with reasonably slow changes in the parameters, the indicators can capture the contemporary stability of the system given by $\lambda$ to a satisfactory degree.

The following settings are considered in this analysis. The linear restoring rate $\lambda$ either follows the decline that is typical for a fold bifurcation in normal form with a linear change in the bifurcation parameter, i.e.,

$$\lambda(t) := \lambda_0 \sqrt{1 - t/T},$$

or it remains constant, i.e.,

$$\lambda(t) \equiv \lambda_0.$$

Here, $T = 14000$ is the time span of the complete experimental setup, and $\lambda_0 \sim \mathcal{U}(0.3, 0.5)$ is a randomly drawn scaling parameter. Note that $\theta$ and $\kappa$ evolve linearly starting from $\theta_0, \kappa_0 \sim \mathcal{U}(0.5, 4)$ and ending at $\theta_T, \kappa_T \sim \mathcal{U}(0.5, 4)$, respectively:

$$\theta(t) := \left(1 - \frac{t}{T}\right)\theta_0 + \frac{t}{T}\theta_T, \qquad \kappa(t) := \left(1 - \frac{t}{T}\right)\kappa_0 + \frac{t}{T}\kappa_T.$$

The samples are generated by discrete-time integration of the continuous-time differential equation in Eq. (2) via the Euler method using time steps $\delta t = 0.1$ after having integrated Eq. (8). In 20 disjoint windows of size $N = 700$ each, we apply the four estimators in question and calculate the Kendall's $\tau$ value for each of these indicator series of size 20. Applying the methods on overlapping windows instead of a disjoint partition does not affect any of the presented results. We draw 5000 random instances of $(\lambda_0, \theta_0, \theta_T, \kappa_0, \kappa_T)$, and for each, we generate one sample time series for a truly decreasing and one for a constant $\lambda(t)$. Based on these findings, we may assess the true- and false-positive rates of the different indicators and, hence, benchmark our newly proposed indicators against existing ones. A visualization of one of these sample instances, along with the relevant parameter thresholds, is given in Fig. 1. The corresponding sample paths in this case clearly show a spurious increase in the observed variance and AC (1), rendering them unsuitable indicators of CSD despite their wide usage. The reason is that the trends in both noise parameters $\theta$ and $\kappa$ have the same effect as decreasing $\lambda$ with respect to these quantities.

The ROC curve is determined by varying the threshold value demanded of the Kendall's $\tau$ to qualify as a significant increase in the respective estimator. If this threshold is high, there will be a high number of false negatives in the decreasing $\lambda$ case. As one gradually lowers the threshold (moving from the bottom left to the top right in Fig. 2), a good indicator will show a more rapid increase in true-positive than in false-positive results, which results in a characteristic arc toward the top-left corner of the plane for a good indicator. A one-dimensional performance metric of the estimator is the area under the ROC curve (AUC), which is a quantity commonly employed for the comparison of CSD indicators [37–40]. The ROC curves, along with the respective AUC values, can be seen in Fig. 2(a).

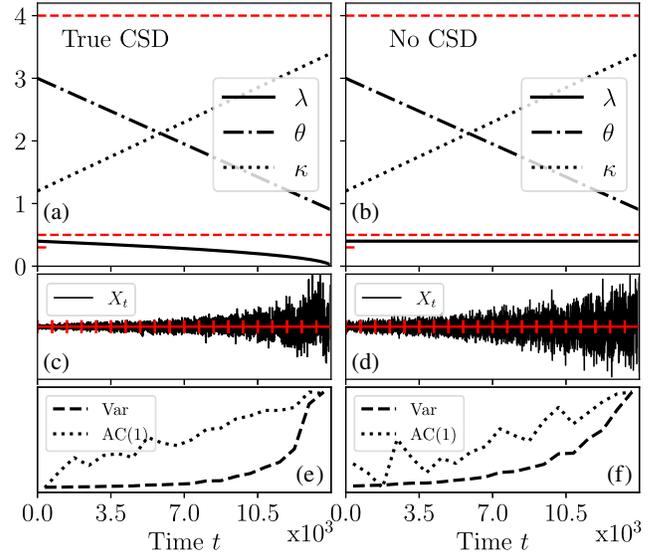

FIG. 1. Example of a randomly generated parameter setting according to the procedure described in the main text. (a),(b) Evolutions of $\lambda(t)$, $\theta(t)$, and $\kappa(t)$. The $\lambda(t)$ corresponding to true critical slowing down is given in panel (a) while that of the null model is given in panel (b). The red dashed lines represent the boundary values of the uniform distributions from which the start and end values of the parameters were drawn. (c),(d) Respective generated sample paths, along with the partition into disjoint windows for the subsequent application of the estimators [results shown in panels (e) and (f)]. In this parameter setting, a spurious indication of CSD is possible, as can be observed in panel (f). The reason is that the trends of $\theta(t)$ and $\kappa(t)$ influence the quantities of variance and AC(1) of the observable $X$ in a way that is indistinguishable from true CSD.

In the context of assessment through Kendall's $\tau$ trends, it bears mentioning that the kind of sensitivity-specificity analysis inherent to the ROC is missing the information on the explicit threshold value along the curve. In order to obtain a complete curve, the threshold value may have to be reduced to $-1$, thus interpreting negative values in Kendall's $\tau$ as positive outcomes. In our plots of the ROC curves, we mark the point along the curve at which the last sensible threshold of zero is crossed. The higher the true-positive rate (TPR) associated with this point, the larger the amount of true-positive classifications that are indeed reasonable. As expected, the false-positive rate (FPR) associated with this point is approximately 50% for all of the indicators since the null model is, by construction, symmetric with respect to the parameter trends (Fig. 2). In an application, one would likely choose a higher classification threshold for the Kendall's $\tau$ value, which would curtail the false-positive rate and increase the significance of positive results in a statistical sense.

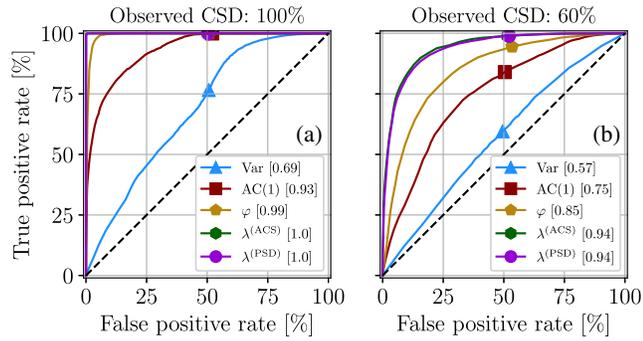

FIG. 2. ROC curves for each of the considered CSD indicators obtained through the procedures laid out in the main text. (a) Comparison performed on the entirety of the decline in $\lambda(t)$ typical for a fold bifurcation. The ROC curves of the indicators using the ACS and PSD concur at nearly perfect discrimination. The AUC is given in the legends. The symbols mark the locations within the curve generation at which the threshold value of the Kendall's $\tau$ is zero. In panel (b), only 60% of the trend and 60% of available data were used. The two novel indicators introduced here, based on the ACS and PSD, respectively, perform best in both comparison settings. The quality of all indicators deteriorates as the amount of available data decreases and the underlying trend in $\lambda(t)$ becomes less pronounced, corresponding to the setting where one is still far away from a bifurcation-induced transition. The individual reasons for this finding are discussed in the main text.

Comparing the AUC values of the five indicators in Fig. 2(a), the conventional marker variance and AC(1) perform worst under our broad model assumption of evolving noise characteristics. However, AC(1) still captures the CSD better than the variance. The estimator $\hat{\varphi}$ and the novel estimation methods via ACS and PSD give the most robust indication of whether CSD is actually taking place in this case.

The sharp dropoff in the linear restoring rate towards the end of the parameter time series is characteristic of fold-type bifurcations. Yet, allowing assessment of CSD up until the bifurcation point can give an unrealistically positive impression of the indicator's skill because, in the immediate proximity to a bifurcation point, noise-induced tipping may become inevitable [41,42]. In applications, the indicator should be able to confidently assess whether CSD is taking place long before the sharp dropoff in $\lambda$ dominates the dynamics. Thus, we are prompted to perform a similar comparison to that above but on data corresponding to the first 60% of the evolution in the linear restoring rate $\lambda$.

The noise parameters $\theta$ and $\kappa$ still evolve according to the same constraints as before but now in only the first 60% of the time frame. The amount of available disjoint windows for the respective estimations also reduces to this fraction accordingly.

We give the ROC curves for this more realistic setting in Fig. 2(b) and further illustrate the comparison between the two situations in Fig. 3. In the more realistic case of having only access to data still relatively far from the bifurcation point, the new CSD indicators introduced here ($\hat{\lambda}^{(\mathrm{ACS})}$ and $\hat{\lambda}^{(\mathrm{PSD})}$) dramatically outperform all other indicators in terms of the ROC curves and the AUC, including $\hat{\varphi}$. It should also be noted that all indicators perform worse in this more realistic comparison setting compared to the case of full access to the data, up to the bifurcation point. For the two conventional indicators variance and AC(1), the reduced performance is due to the relative sizes of the parameter trends. Lacking the sharp decline in $\lambda$ toward the bifurcation, the changes in the other two parameters can more easily overwhelm the effects of a changing $\lambda$. For the indicators proposed here, using the estimators $\hat{\lambda}^{(\mathrm{ACS})}$ and $\hat{\lambda}^{(\mathrm{PSD})}$, the reason for the increase in faulty results is not rooted in spurious or masking effects themselves. Instead, a higher uncertainty associated with the estimations leads to more noise in the indicator trends. This uncertainty stems from the fact that the methods respond more sensitively to the fast parameter changes in $\theta$ and $\kappa$. A larger amount of longer windows of data would work against this statistical effect, yet in applications, the amount of data available is often limited and of the order of magnitude discussed in this section. Thus, comparing the techniques boils down to a trade-off between exposure to spurious indication or masking and potential misestimation due to a lower signal-to-noise ratio. Nevertheless, our results show that the two indicators proposed in this work should always be preferred over the conventional variance and AC(1).

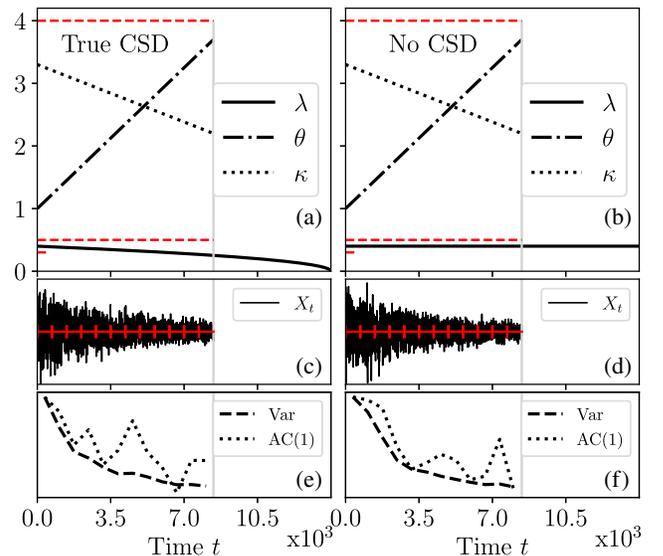

FIG. 3. Example for a randomly generated parameter setting according to the intensified comparison procedure described in the main text (panels are equivalent to those in Fig. 1), but note the inverted trends of $\theta(t)$ and $\kappa(t)$. In this parameter setting, masking of CSD is possible. The reason is that the trends of $\theta(t)$ and $\kappa(t)$ counteract the effect of a decreasing $\lambda(t)$ in the theoretical stationary quantities of variance and AC(1) of the observable $X$.

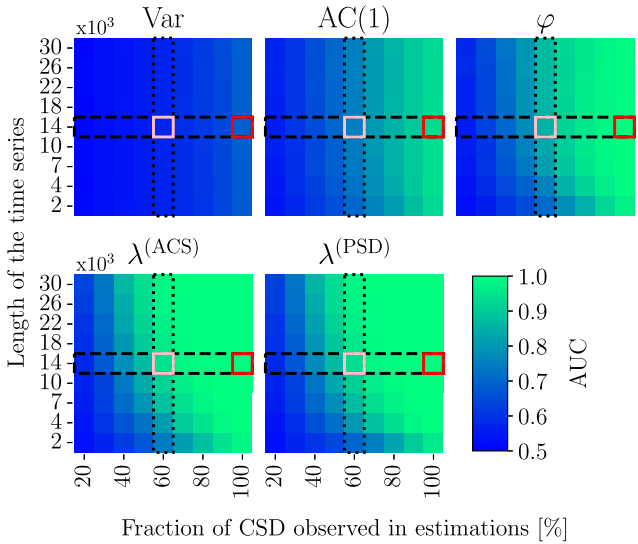

FIG. 4. AUC values obtained from ROC curves of each of the discussed CSD indicators in different observational settings. The parameter evolutions of $\lambda(t)$, $\theta(t)$, and $\kappa(t)$ are generated as outlined in the main text. The length of the time series indicated on the $y$ axis of each panel is divided into 20 disjoint windows. The percentage on the $x$ axis determines how many of these windows are used to assess CSD in the synthetically generated data. Highlighted in pink and red are the values corresponding to the ROC curves presented in Figs. 2(a) and 2(b), respectively. The black cross sections of these heat maps are shown in Figs. 5(a) (dashed) and 5(b) (dotted), respectively.

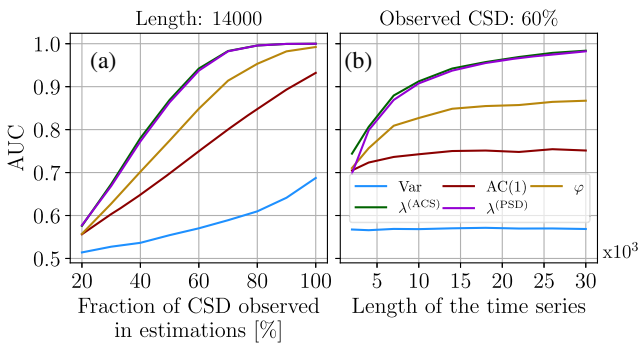

FIG. 5. One-dimensional cross sections of the heat maps in Fig. 4. In panel (a), the length of the evaluated time series is kept constant at $T = 14000$ while the fraction of the time series before tipping that is used for CSD estimation with the different indicators is varied. In panel (b), this fraction is kept constant at 60% while the length of the evaluated time series is varied.

We further illustrate this trade-off by performing the above analysis for an ensemble of different configurations, varying the time-series length and the percentage of the time series used to detect CSD. Figure 4 shows the AUC values for each CSD indicator under these varying conditions, and Fig. 5 shows cross sections along fixed choices of the total time-series length and the observed fraction. It is apparent that the quality of assessment for the conventional methods almost exclusively depends on the fraction of CSD under observation and not the amount of data points. The novel methods designed for the continuous-time red noise case can detect CSD at much earlier points in time, i.e., after very few observed windows. Even though the discrete-time red noise method via $\hat{\varphi}$ performs better than the two conventional indicators, the novel methods still clearly outperform it.

### B. Analyzing the desertification of the Western Sahara

In the following, we show how the methods introduced above can be used to discriminate between different physical candidate mechanisms leading to real-world abrupt transitions, focusing on the example of the abrupt desertification of the Western Sahara some 6 thousand years ago. Hopcroft and Valdes [29] investigated the retreat of Western Sahara vegetation during the mid-Holocene epoch and, in particular, whether climate models support the view of an abrupt retreat possibly caused by bifurcation dynamics. In the Western Sahara region, the contemporary climate is that of an arid, hot desert. Paleoclimate evidence suggests that during the late Pleistocene and early to mid-Holocene epochs, until about 6 thousand years ago, there was abundant savanna-type vegetation present in the same region [43]. The driving external variable responsible for this change was the orbital forcing, which affects the Northern Hemisphere summer insolation [44]. Before 6 thousand years ago, the increased summer insolation in the Western Sahara facilitated the green Sahara via the following feedback mechanism: The vegetation in the region had a lower reflectivity than the desert and hence absorbed more solar energy, which can fuel convective systems and even cause a northward extension of the West African summer monsoon system [45]. Changes in cloud cover and evapotranspiration must also be considered [45,46].

There has been a debate about whether the available paleoclimate data allow for the characterization of the Western Sahara vegetation system as a tipping element [47–49], in the strict sense of exhibiting bifurcation dynamics that can lead to critical transitions between alternative states. Even though the aforementioned conceptual "Charney" model of the above-described feedback is plausible, there may be more complex and spatially constrained dynamics at play. A suitable consistency check of the hypothesis of positive feedback mechanisms driving the transition is to investigate the time series for indications of critical slowing down. Hopcroft and Valdes [29,50] performed a preliminary analysis of the variance in vegetation coverage in the advent of the transition and found a clear increase. This finding may be interpreted as an indication of critical slowing down if the underlying model assumptions on the noise are valid. These include the rather strict assumption of the disturbances inflicted on the system being well represented by stationary white noise and excludes nonstationary temporal correlations. Such correlations can be found in atmospheric observables, which are relevant to the dynamics of vegetation systems. We will analyze the time-series data obtained from the model configuration of Hopcroft and Valdes [29] with respect to critical slowing down using our novel estimation methods for quantifying system stability. The underlying premise for the applicability of the red noise model is to assume that disturbances in precipitation drive the vegetation dynamics. More specifically, the following linearized conceptual model of the Western Sahara vegetation system warrants the application of our novel methods to simulated data from the complex general circulation model of Ref. [29] (see SM [30] for a more detailed motivation [30]):

$$dV_t = -\lambda(t)\big(V_t - \overline{V}(t)\big)dt + \kappa_V(t)\big(P_t - \bar{P}(t)\big)dt, \quad (9a)$$

$$dP_t = -\theta(t)\big(P_t - \bar{P}(t)\big)dt + \kappa_P(t)dW_t, \quad (9b)$$

where $\bar{V}(t)$ denotes the equilibrium of $V(t)$. The overall negative feedback strength $-\lambda(t)$ acts on the yearly averaged spatially extended vegetation $V(t)$, measured as the fraction of ground covered by certain plant functional types. Furthermore, the rate of change of $V$ away from equilibrium $\bar{V}(t)$ is assumed to be proportional to deviations of the precipitation $P$ from its contemporary equilibrium $\bar{P}(t)$. These deviations are, in turn, modeled as an Ornstein-Uhlenbeck process with correlation parameter $\theta(t)$. The equilibrium dynamics of the vegetation thus follow the linearized model driven by continuous-time red noise [cf. Eq. (1)] introduced and analyzed in the previous sections. The conceptual model introduced in Eqs. (9a) and (9b), though not further investigated numerically, serves as the argumentative basis for the application of our methods to the time-series data of Ref. [29]. A decrease of the feedback parameter $\lambda(t)$ can, in this context, be interpreted as a weakening of the stability of the vegetation system in the general circulation model. Such indications on the basis of CSD, were they to be found, would imply that positive feedback is gaining in strength while negative feedback weakens.

On yearly averaged time-series data of vegetation and precipitation obtained from the climate model in Ref. [29], we first determine $\bar{V}(t)$ and $\bar{P}(t)$ by applying a Gaussian filter. The analysis is performed up until the observed tipping point, which is defined to be the point of highest negative curvature in $\bar{V}(t)$. The stability estimators based on the ACS and PSD are each employed on the time-series data of $V_t - \bar{V}(t)$ and $P_t - \bar{P}(t)$. In the first case, $\lambda(t)$ and $\theta(t)$ are inferred from the vegetation data, while in the second case, $\theta(t)$ is inferred from the precipitation data. A consistency check of the presuppositions made in our conceptual model defined by Eqs. (9a) and (9b) can be performed by comparing the two estimations of $\theta(t)$. Figure 6 shows the results for the data of one specific simulation grid cell at 25 °N and 3.75 °W. Analogous analyses with similar results for other grid cells can be found in the SM [30]. A decrease in $\lambda(t)$ can clearly be observed in all of these applications. Thus, the increase in variance in the advent of the transition observed by Hopcroft and Valdes [29,50], using the methodology introduced here, can be attributed to an actual decrease in system stability. The results from $\theta(t)$ stemming from the two time series match qualitatively, encouraging the proposed conceptual linear model as a good representation of the equilibrium dynamics.

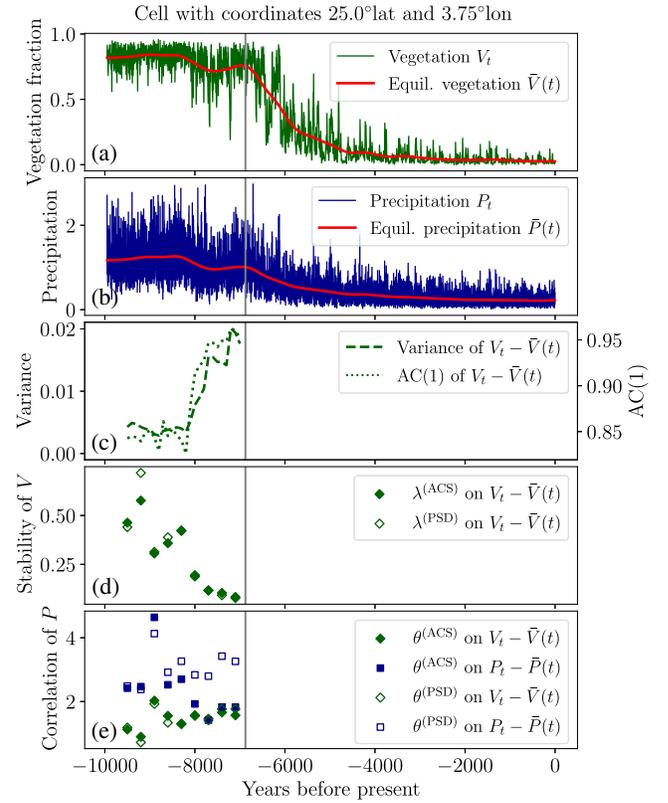

FIG. 6. Analysis of simulations from a complex climate model for the African Humid Period obtained from Ref. [29] for the model cell located at 25° latitude and 3.75° longitude. Panels (a) and (b) show the time-series data and its contemporary mean obtained via a Gaussian filter. The tipping point is marked in gray. In panel (c), the conventional EWS of variance and AC(1) suggest CSD. The estimation of system stability $\lambda$ via the ACS and PSD of $V$ in panel (d) supports this indication consistently and rules out the counter-hypothesis, namely, effects of nonstationary red noise as the main cause. We emphasize that this discrimination could not be performed based on variance and AC(1). In panel (e), the two estimations of the correlation parameter $\theta$ of the precipitation based on data from $V$ and $P$, respectively, can be compared. While they differ at times by a factor of about 2, they match qualitatively in order of magnitude. Their individual trends are approximately flat, indicating no substantial change in the correlation characteristics over the observed time span.

## IV. DISCUSSION

The estimators for variance and AC(1) commonly employed as indicators for CSD easily lead to a false assessment when aspects of the driving noise cannot be assumed to be constant. In the case of the general red noise model, we have discussed this finding on the basis of theoretical considerations and demonstrated it on sample data.

The two, new CSD indicators we introduced here, designed to be sensitive to changes in the correlation characteristics of the red noise, perform substantially better across a broad range of parameter configurations as measured by the receiver-operating characteristic. However, their performance still depends on the length of the given time-series data, as seen in Figs. 4 and 5. In effectively every configuration of the size and number of observed windows given there, the two novel methods outperform other existing methods of detecting CSD, including methods designed for discrete-time red noise.

Many questions in the context of potentially bifurcation-induced tipping in applications may be more robustly assessed with these new methods. We presented the example of the desertification of the Green Sahara. Applying our methods to paleoclimate model data reveals that this archetype of abrupt climate change is indeed associated with a weakening of negative feedback, in the underlying physical system. This finding has crucial implications for the existence of similar instabilities in the current climate and the potential crossing of tipping points in the context of ongoing anthropogenic global warming.

We stress that, in general, the rather specific red noise model need not be a good fit for observed time series without first performing an adequate analysis. In most cases, it will rely on a physical understanding of the underlying dynamics. In order to apply our methods to data concerning the desertification of the Western Sahara, we have posited such a conceptual model and performed an analysis of model consistency and system stability based on the ACS and PSD of the available model data. The results allow for the attribution of the previously observed increase in variance before the transition to a destabilization of the system measured by its linear restoring rate. In the absence of such confirmation, changes in the driving noise of the system could not be excluded as the main cause for observed critical slowing down.

The GitHub repository RedNoiseEstimatorComparison [51] can be used to access the code generating all figures in this manuscript. Therein, numerical implementations of all discussed methods are included.

## ACKNOWLEDGMENTS

We thank P. Hopcroft and P. Valdes for helpful discussions and for providing the climate model simulations. This work is ClimTip contribution #5; the ClimTip project has received funding from the European Union's Horizon Europe research and innovation programme under Grant Agreement No. 101137601. N. B. acknowledges further funding from the Volkswagen Stiftung and the European Union's Horizon 2020 research and innovation programme under the Marie Sklodowska-Curie Grant Agreement No. 956170, funded by the European Union.

Views and opinions expressed are those of the author(s) only and do not necessarily reflect those of the European Union or the European Climate, Infrastructure and Environment Executive Agency (CINEA). Neither the European Union nor the granting authority can be held responsible for views and opinions in this work.